\documentclass[a4paper,leqno,11pt,final]{article}
\usepackage{latexsym}
\usepackage{amsthm,amsmath,amssymb}

\usepackage[dvips]{graphicx, color}

\newtheorem{thm}{Theorem}

\newtheorem{prop}[thm]{Proposition}

\newtheorem{conj}[thm]{Conjecture}
\newtheorem{lem}[thm]{Lemma} 

\newcommand{\Aut}{\mathrm{Aut}}

\newcommand{\z}{\widehat\zeta_E}

\begin{document}
\title{{\bf\Large{Zeta Functions for Elliptic Curves}}\\
{\bf\large{I. Counting Bundles}}}
\author{\bf Lin WENG}
\date{\bf }
\maketitle
\noindent
{\bf Abstract:} {\footnotesize To count bundles on curves, we study zetas of elliptic curves and their zeros. There are two types, i.e., the pure non-abelian zetas defined using moduli spaces of semi-stable bundles, and the group zetas defined for special linear groups. In lower ranks, we show that these two types of zetas coincide and satisfy the Riemann Hypothesis. For general cases, exposed is an intrinsic relation on automorphism groups of semi-stable bundles over elliptic curves, the so-called counting miracle.
All this, together with Harder-Narasimhan, Desale-Ramanan and Zagier's result, gives an effective way to
count semi-stable bundles on elliptic curves not only in terms of automorphism groups but more essentially in terms of their $h^0$'s.
Distributions of  zeros of high rank zetas are also discussed.}
\vskip 1.50cm
$$\mathrm{\bf CONTENTS}$$
\vskip 0.30cm
\noindent
\S1. {\bf High Rank Zetas for Elliptic Curves}

1.1 Fat Moduli Spaces

1.2 Definition

1.3 Rank Two

1.4 Rank Three

\vskip 0.30cm
\noindent
\S2. {\bf Zetas of Elliptic Curves Associated to $SL_n$}

2.1 Definition

2.2 $SL_2$ 

2.3 $SL_3$

2.4 Distribution of zeros

\vskip 0.30cm
\noindent
\S3. {\bf Counting Bundles}

\vskip 0.30cm
\noindent
\S4. {\bf Distribution of Zeros}

4.1 The Riemann Hypothesis

4.2 Distribution of zeros
\section{High Rank Zetas for Elliptic Curves} 
\subsection{Fat Moduli Spaces}

Let $X$ be an irreducible, reduced, regular projective curve of genus $g$ defined over $\mathbb F_q$. Denote by $\mathcal M_{X,r}(d)$ the moduli space of rank $r$ semi-stable bundles of degree $d$ consisting of the Seshadri Jordan-H\"older equivalences of $\mathbb F_q$-rational semi-stable bundles. For our purpose, we consider $\mathcal M_{X,r}(d)$ in the sense of the {\it fat} moduli, meaning that ordinary moduli spaces equipped with the additional structure at Seshadri class $[\mathcal E]$
defined by the collection of semi-stable bundles in $[\mathcal E]$, namely, the set $\big\{\mathcal E:\mathcal E\in [\mathcal E]\big\}$ is adding at the point $[\mathcal E]$. $\mathcal M_{X,r}(d)$ equipped with such a structure is called a fat moduli space and denoted as ${\bf M}_{X,r}(d)$. 

A natural question is to count these 
$\mathbb F_q$-rational semi-stable bundles $\mathcal E$. For this purpose, two invariants, namely, the automorphism group $\Aut(X,\mathcal E)$ and its global sections $h^0(X,\mathcal E)$ can be naturally used.
This then leads to the refined Brill-Noether loci $$W_{X,r}^i(d):=\Big\{[\mathcal E]\in {\bf M}_{X,r}(d):
\min_{\mathcal E\in[\mathcal E]}:h^0(X,\mathcal E)\geq i\Big\}$$ and $$[\mathcal E]^j:=\{\mathcal E\in[\mathcal E]: \dim_{\mathbb F_q}\Aut{\mathcal E}\geq j\}.$$ Recall that there exist natural
isomorphisms $${\bf M}_{X,r}(d)\to {\bf M}_{X,r}(d+rm),\qquad
\mathcal E\mapsto A^m\otimes\mathcal E$$ and $${\bf M}_{X,r}(d)\to {\bf M}_{X,r}(-d+r(2g-2)),\qquad
\mathcal E\mapsto K_X\otimes\mathcal E^\vee,$$ where $A$ is an Artin line bundle of degree one and $K_X$ denotes the dualizing bundle of $X/\mathbb F_q$. So we only need to count ${\bf M}_{X,r}(d_0)$ for $d_0= 0, 1,\dots, r(g-1).$ Accordingly, we set
$$\alpha_{X,r}(d):=\sum_{\mathcal E\in {\bf M}_{X,r}(d)}\frac{q^{h^0(X,\mathcal E)}-1}{\#\Aut(\mathcal E)},\qquad \beta_{X,r}(d):=\sum_{\mathcal E\in {\bf M}_{X,r}(d)}\frac{1}{\#\Aut(\mathcal E)}$$ with $\beta$ a classical invariant ([HN]).

So, to count bundles, the problem becomes how to control $\alpha_{X,r}(d_0)$'s with $d_0$ ranging as above, and 
$\beta_{X,r}(q)$ with $q=0,1,\dots r-1$. For $\alpha$, two general principles can be used for counting semi-stable bundles, namely, the vanishing theorem claiming that
$h^1(X,\mathcal E)=0$ if $d(\mathrm E)\geq r(2g-2)+1$ and the Clifford lemma claiming that
$h^0(X,\mathcal E)\leq r+\frac{d}{2}$ if $0\leq \mu(\mathcal E)\leq 2g-2$. But this is merely the starting point.
By contrasting, the invariant $\beta$ can be understood, thanks to the high profile works of Harder-Narasimhan ([HN]),
Desale-Ramanan ([DR]) and Zagier ([Z]). 

To state it, let $$\zeta_X(s):=\frac{\prod_{i=1}^{2g}(1-\omega_iq^{-s})}{(1-q^{-s})(1-qq^{-s})}$$ be the Artin zeta function of $X/\mathbb F_q$,
 $$v_r(\mathbb F_q):=v_r:=\frac{\prod_{i=1}^{2g}(1-\omega_i)}{q-1}q^{(r^2-1)(g-1)}\zeta_X(2)\cdots\zeta_X(r)$$ and
 $$c_{r,d}(q):=c_{r,d}:=\prod_{i=1}^{s-1}\frac{q^{(n_i+n_{i+1})\{n_1+\cdots+n_i)d/n\}}}{1-q^{n_i+n_{i+1}}}.$$  With above, then the above works mentioned can be strengthen as follows:
  
 \begin{thm} ([Z, Thm 2]) For any pair $(r,d)$, we have
 $$\beta_{X,r}(d)=\sum_{n_1,\dots,n_s>0,\sum n_i=r}q^{(g-1)\sum_{i<j}n_in_j} c_{r,d}(q)\prod_{i=1}^s v_{n_i}(\mathbb F_q).$$
 \end{thm}
 
 \subsection{Definition}

Let $E$ be a regular, integral projective elliptic curve defined over $\mathbb F_q$.
Define {\it rank r pure zeta function} for $E/\mathbb F_q$ by
$$\begin{aligned}\widehat\zeta_{E,r}(s):=\zeta_{E,r}(s):=&\sum_{m=0}^\infty
\sum_{V\in {\bf  M}_{E,r}(d), d=rm}\frac{q^{h^0(C,V)}-1}{\#\mathrm{Aut}(V)}\cdot (q^{-s})^{d(V)}.\end{aligned}$$ 
Then, by the vanishing theorem for semi-stable bundles, 
$$\begin{aligned}\zeta_{E,r}(s)=&\sum_{V\in {\bf  M}_{E,r}(0)}\frac{q^{h^0(C,V)}-1}{\#\mathrm{Aut}(V)}+\sum_{m=1}^\infty
\sum_{V\in {\bf  M}_{E,r}(rm)}\frac{q^{h^0(C,V)}-1}{\#\mathrm{Aut}(V)}\cdot (q^{-s})^{d(V)}\\
=&\alpha_{E,r}(0)+\sum_{m=1}^\infty
\sum_{V\in {\bf  M}_{E,r}(rm)}\frac{q^{rm}-1}{\#\mathrm{Aut}(V)}\cdot (q^{-s})^{rm}\\
=&\alpha_{E,r}(0)+\beta_{E,r}(0)\sum_{m=1}^\infty
(q^{rm}-1)\cdot (q^{-s})^{rm}\\
=&\alpha_{E,r}(0)+\beta_{X,r}(0)\cdot \Big(\frac{(qt)^r}{1-(qt)^r}-\frac{t^r}{1-t^r}\Big)
\end{aligned}$$
Consequently,
$$\widehat Z_{E,r}(\frac{1}{qt})=\widehat Z_{E,r}(t)
=\alpha_{X,r}\big(0\big)+\beta_{E,r}(0)\cdot \frac{(Q-1)T}{(1-T)(1-QT)}.$$
Here, 
$T:=t^r,\, Q=q^r$, $Z_{E,r}(t)=\zeta_{E,r}(s)$, and $\widehat Z_{E,r}(t)=\widehat\zeta_{E,r}(s)$. This then completes the proof of the following

\begin{thm} (i) $\zeta_{E,1}(s)=\zeta_E(s)$, the Artin zeta function for $E/\mathbb F_q$;

\noindent
(ii) ({\bf Rationality}) There exists a degree $2$ polynomial $P_{E,r}(T)\in\mathbb Q[T]$ of $T$
such that
$$Z_{E,r}(t)=\frac{P_{E,r}(T)}{(1-T)(1-QT)}\quad\mathrm{with}\quad T=t^r,\ Q=q^r;$$
\noindent
(iii) ({\bf Functional equation}) $$\widehat Z_{E,r}(\frac{1}{qt})=\widehat Z_{E,r}(t),$$
\end{thm}

\noindent
{\it Remark.} The pure zeta here, a new genuine one, is quite different from the zeta introduced in [W1]. The reason for the purity is that the zeta in [W1] does not satisfy the Riemann Hypothesis.

\subsection{Rank Two}
For rank two pure zeta, it suffices to calculate 
$\alpha_{E,2}(0)$ and $\beta_{E,2}(0)$. For $\beta$,
 by Harder-Narasimhan, Desale-Ramanan and Zagier's formula, i.e., Thm 1, 
 $$\beta_{E,2}(0)=\frac{N}{q-1}\Big(1+\frac{N}{q^2-1}\Big).$$
 Here, as usual,  $N$ denotes the number of $\mathbb F_q$-rational points of $E$.
On the other hand, by the classification of Atiyah ([A]), over $\overline{\mathbb F_q}$, the graded bundle associated to a Jordan-H\"older filtration of a semi-stable bundle $V\otimes_{\mathbb F_q}\overline{\mathbb F_q}$ is of the form $\mathrm{Gr}(V\otimes_{\mathbb F_q}\overline{\mathbb F_q})=L_1\oplus L_2$ with $L_i$ degree zero line bundles, which may not really defined over $\mathbb F_q$. Consequently, for $\mathbb F_q$-rational semi-stable bundles $V$ of rank two, $h^0(E,V)\not=0$ if and only if $V=\mathcal O_E\oplus L$ or $V=I_2$ with $L$ a $\mathbb F_q$-rational line bundle of degree 0 and $I_2$ the only non-trivial extension of $\mathcal O_E$ by itself. Thus $\alpha_{E,2}(0)$ is given by $$\begin{aligned}&\Big(\frac{q^{h^0(E,\mathcal O_E\oplus \mathcal O_E)}-1}{\#\Aut(\mathcal O_E\oplus \mathcal O_E)}+\frac{q^{h^0(E,I_2)}-1}{\#\Aut (I_2)}\Big)+\sum_{L\in\mathrm{Pic}^0(E), L\not=\mathcal O_E}\frac{q^{h^0(E,\mathcal O_E\oplus L)}-1}{\#\Aut(\mathcal O_E\oplus L)}\\
=&\Big(\frac{q^2-1}{(q^2-1)(q^2-q)}+\frac{q-1}{(q-1)q}\Big)+(N-1)\frac{q-1}{(q-1)^2}=\frac{N}{q-1}.\end{aligned}$$ Thus, we have the following

\begin{prop} (i) $\displaystyle{\alpha_{E,2}(0)=\beta_{E,1}(0)=\frac{N}{q-1}}$;
\vskip 0.20cm
\noindent
(ii) $\displaystyle{Z_{E,2}(t)=\alpha_{E,2}(0)\cdot\frac{1+(N-2)T+QT^2}{(1-T)(1-QT)}}.$
\end{prop}
The Riemann Hypothesis holds since $$\Delta_2=(N-2)^2-4Q=(N-2-2q)(N-2+2q)<0$$ using Hasse's theorem for the Riemann Hypothesis  of elliptic curves, namely $N\leq 2\sqrt q.$ That is to say, we have proved the following
\begin{thm}
({\bf Riemann Hypothesis}$_2$) $$\widehat \zeta_{E,2}(s)=0\qquad\Rightarrow\qquad \mathrm{Re}(s)=\frac{1}{2}.$$
\end{thm}

\subsection{Rank Three}
To understand rank 3 zeta for elliptic curves, we need to calculate $\alpha_{E,3}(0)$. 
Recall that for semi-stable vector bundles $V$ over elliptic curves, by Atiyah ([A]), the Jordan-H\"older graded bundles $G(V)$ are of the form $L_1\oplus L_2\oplus L_3$ with $L_i\in\mathrm{Pic}^0(E\otimes\overline{\mathbb F_q})$. In particular, for semi-stable bundles with non-trivial contribution to $\alpha_{E,3}(0)$, at least one of the $L_i$'s should be $\mathcal O_E$. So three types:

(1) $L_1=L_2=L_3=\mathcal O_E$;

(2) $L_1=L_2=\mathcal O_E\not=L_3$;

(3)  $L_1=\mathcal O_E\not=L_i,\, i=2,\,3.$

Accordingly, write $$\displaystyle{\alpha_{E,3}(0)=\alpha_{E,3}^1(0)+\alpha_{E,3}^2(0)+\alpha_{E,3}^3(0),}$$ with $$\alpha_{E,3}^i(0)=\sum_{\substack{V\in{\bf M}_{E,3}(0),\\ V:\ \mathrm{type}\ i}}\frac{q^{h^0(X,V)}-1}{\#\Aut(V)},\qquad 1=1,2,3.$$

Clearly, for types (1) and (2), all $L_i$'s are $\mathbb F_q$-rational. However, for type (3), it may well be possible that $L_2$ and $L_3$
are not, even $V$ itself is $\mathbb F_q$-rational. For general curves, this proves to be an essential difficulty in understanding rank 3 zetas.
However, for elliptic curves, due to the degree constraint, we find a nice way to overcome it.

The idea is to use Thm 1. To be more precise, for type (3), we have
$$\begin{aligned}\alpha_{E,3}^3(0)=&\sum_{\substack{V\in{\bf M}_{E,3}(0),\\ G(V)=\mathcal O_E\oplus L_2\oplus L_3\\
L_2\not=\mathcal O_E\not=L_3}}\frac{q^{1}-1}{\#\Aut(V)}\\
=&\sum_{\substack{W\in{\bf M}_{E,2}(0),\\ G(W)= L_2\oplus L_3\\
L_2\not=\mathcal O_E\not=L_3}}\frac{q-1}{\#\Aut(W)\cdot(q-1)}
=\sum_{\substack{W\in{\bf M}_{E,2}(0),\\ G(W)= L_2\oplus L_3\\
L_2\not=\mathcal O_E\not=L_3}}\frac{1}{\#\Aut(W)},\end{aligned}$$
since there is no non-zero morphisms between $\mathcal O_E$ and $L_i$'s.
But $$\begin{aligned}&\Big\{W\in{\bf M}_{E,2}(0): G(W)= L_2\oplus L_3,
L_2\not=\mathcal O_E\not=L_3\Big\}\\
&={\bf M}_{E,2}(0)\Big\backslash\Big\{W\in{\bf M}_{E,2}(0):G(W)=\mathcal O_E\oplus L \Big\}.\end{aligned}$$
Consequently,
$$\begin{aligned}\alpha_{E,3}^3(0)=&\beta_{E,2}(0)-
\sum_{\substack{W\in{\bf M}_{E,2}(0),\\ G(W)= \mathcal O_E\oplus L}}\frac{1}{\#\Aut(W)}\\
=&\beta_{E,2}(0)-
\sum_{\substack{W\in{\bf M}_{E,2}(0),\\ G(W)= \mathcal O_E\oplus L,\, L\not=\mathcal O_E}}\frac{1}{\#\Aut(W)}-\Big(\frac{1}{\#\Aut(\mathcal O_E\oplus \mathcal O_E)}+\frac{1}{\#\Aut(I_2)}\Big)\\
=&\beta_{E,2}(0)-
\sum_{\substack{W= \mathcal O_E\oplus L,\, L/\mathbb F_q\not=\mathcal O_E}}\frac{1}{\#\Aut(W)}-\Big(\frac{1}{(q^2-1)(q^2-q)}+\frac{1}{(q-1)q}\Big)\\
=&\beta_{E,2}(0)-\frac{N-1}{(q-1)^2}-\frac{q}{(q^2-1)(q-1)}\\
=&\frac{N^2+(q^2-q-2)N+1}{(q^2-1)(q-1)}.\end{aligned}$$
Now by Thm 1, we have $$\begin{aligned}\beta_{E,2}(0)=&\frac{N}{q-1}\zeta_E(2)+\frac{1}{1-q^2}\Big(\frac{N}{q-1}\Big)^2\\
=&
\frac{N}{q-1}\Big(1+\frac{N}{q^2-1}\Big).\end{aligned}$$

For type (2), we have $$\begin{aligned}\alpha_{E,3}^2(0)=&\sum_{L/\mathbb F_q\not=\mathcal O_E}\Big(\frac{q^2-1}{\#\Aut(\mathcal O_E\oplus \mathcal O_E\oplus L)}+\frac{q-1}{\#\Aut(I_2\oplus L)}\Big)\\
=&(N-1)\Big(\frac{q^2-1}{(q^2-1)(q^2-q)(q-1)}+\frac{q-1}{(q^2-q)(q-1)}\Big)\\
=&\frac{q}{(q^2-1)(q-1)}.\end{aligned}$$

Finally for type (1), we have 
$$\begin{aligned}\alpha_{E,3}^1(0)=&\frac{q^3-1}{\#\Aut(\mathcal O_E\oplus \mathcal O_E\oplus \mathcal O_E)}+\frac{q^2-1}{\#\Aut(I_2\oplus \mathcal O_E)}+\frac{q-1}{\#\Aut(I_3)}\\
=&\frac{q^3-1}{(q^3-1)(q^3-q)(q^3-q^2)}+\frac{q^2-1}{(q-1)^2q^3}+\frac{q-1}{(q-1)q^2}\\
=&\frac{N-1}{(q-1)^2}.\end{aligned}$$

Put all this together, we have

\begin{prop} (i)
$\alpha_{E,3}(0)=\beta_{E,2}(0)=\frac{N}{q-1}\cdot\Big(1+\frac{N}{q^2-1}\Big);$

\noindent
(ii) $\beta_{E,3}(0)=\frac{N}{q-1}\cdot\Big[1+\frac{q+2}{q^3-1}N+\frac{N^2}{(q^3-1)(q^2-1)}
\Big];$

\noindent
(iii) $\widehat Z_{E,3}(t)=\frac{N}{q-1}\cdot\frac{1}{(1-T)(1-q^3T)}$\\
$\hskip 3.0cm\times\Big[\Big(1+\frac{N}{q^2-1}\Big)(1+q^3T^2)+\Big(-2+\frac{(2q-3)N}{q-1}+\frac{N^2}{q^2-1}\Big)\cdot T\Big].$
\end{prop}

\noindent
{\it Proof.} For our convenience, set $\z(1)=\frac{N}{q-1}$. Then for (ii), by Thm 1 and $\z(3)=\z(1)\cdot\Big(1+\frac{q^2N}{(q^3-1)(q^2-1)}\Big)$, we have
$$\begin{aligned}\beta_{E,3}(0)=&\z(1)\z(2)\z(3)-\frac{2}{1-q^3}\z(1)\z(1)\z(2)\\
&\qquad+\frac{1}{(1-q^2)(1-q^2)}\z(1)\z(1)\z(1)\\
=&\z(1)\Big[1+\frac{q+2}{q^3-1}N+\frac{N^2}{(q^3-1)(q^2-1)}
\Big].\end{aligned}$$

For (iii), from the expression of $\zeta_{E,r}(s)$ obtained in Thm. 2, 
$$\begin{aligned}\widehat Z_{E,3}(t)
=\frac{N}{q-1}\qquad&\\
\times\Big[\Big(1+\frac{N}{q^2-1}\Big)&+\Big(1+\frac{q+2}{q^3-1}N+\frac{N^2}{(q^3-1)(q^2-1)}
\Big)\frac{(q^3-1)T}{(1-T)(1-q^3T)}\Big]\\
=\frac{\z(1)}{(1-T)(1-q^3T)}&\\
\qquad\times\Big[\Big(1+\frac{N}{q^2-1}\Big)&(1+q^3T^2)+\Big(-2+\frac{(2q-3)N}{q-1}+\frac{N^2}{q^2-1}\Big)T\Big]\\\end{aligned}$$

\section{Zetas of Elliptic Curves associated to $SL_n$}
\subsection{Definition}
For $G=SL_n$ with $B$ the standard Borel subgroup consisting of upper triangular matrices, let $T$ be the associated torus consisting of diagonal matrices. Then the root system $\Phi$ associated to $T$ can be realized 
as $$\Phi^+=\{e_i-e_j:1\leq i<j\leq n\}$$ with $\{e_i\}_{i=1}^n$ the standard orthogonal basis of the Euclidean space $V=\mathbb R^n$. Being type $A_{n-1}$, its simple roots are given by
$$\Delta:=\{\alpha_i:=e_i-e_{i+1}:i=1,2,\dots,n-1\},$$
the so-called Weyl vector is simply $$\rho:=\frac{1}{2}\Big((n-1)e_1+(n-3)e_2+\cdots-(n-3)e_{n-1}-(n-1)e_n\Big),$$ and the Weyl group $W$ may be identified with the permutation group $S_n$ via the action on the subindex of $e_i$'s. Introduce the corresponding fundamental weights $\lambda_j$'s via
$$\langle\lambda_i,\alpha^\vee_j\rangle=\delta_{ij},\qquad\forall\alpha_j\in \Delta.$$
For each $w\in W$, set $\Phi_w:=\Phi^+\cap w^{-1}\Phi^-$.
For $\lambda\in V_{\mathbb C}$,
introduce then the period of $SL_n$ by
$$\omega^{SL_n}_E(\lambda):=\sum_{w\in W}\frac{1}{\prod_{\alpha\in\Delta}(1-q^{-\langle w\lambda-\rho,\alpha^\vee\rangle})}\prod_{\alpha\in \Phi_w}
\frac{\z(\langle\lambda,\alpha^\vee\rangle)}{\z(\langle\lambda,\alpha^\vee\rangle+1)}.$$

Corresponding to $\alpha_P=\alpha_{n-1}$, let $$P=P_{n-1,1}=\Big\{\begin{pmatrix}A&B\\ 0&D\end{pmatrix}:A\in GL_{n-1},D\in GL_1\Big\}$$
be the standard parabolic subgroup of $SL_n$ attacted to the partition $n=(n-1)+1$. 
Write $$\lambda:=\rho+\sum_{j=1}^{n-1}s_j\lambda_j$$ and set $s:=s_{n-1}$. Then 
 we introduce the {\it period for $(SL_n,P)$}
as a one variable function defined by
$$\omega_E^{SL_n/P}(s):=\mathrm{Res}_{s_1=0}\mathrm{Res}_{s_2=0}\cdots
\mathrm{Res}_{s_{n-2}=0}\,\omega^{SL_n}_E(\lambda).$$
This period consists many terms, each of which is a product of certain rational factors of $q^{-s}$ and Atrin zetas. Clear up all zeta factors in the denominators of all terms! The resulting function
is then defined to be the {\it zeta function $\z^{SL_n}(s)$ of $E$ associated to $SL_n$}.

\begin{thm}([W4]) (i) $\z^{SL_n}(s)$ is a well-defined meromorphic function on the whole $s$-plane;

\noindent
(ii) ({\bf  Functional Equation})
$$\z^{SL_n}(-n-s)=\z^{SL_n}(s).$$
\end{thm}

This group theoretic zeta function is expected to play a central role in counting bundles. In fact, we have the following

\begin{conj} (i) ({\bf The Riemann Hypothesis})
$$\z^{SL_n}(s)=0\qquad\Rightarrow\qquad \mathrm{Re}(s)=-\frac{n}{2}.$$

\noindent
(ii) ({\bf Uniformity})
Up to a rational function factor of $q$, $$\z^{SL_n}(-ns)=\widehat\zeta_{E,r}(s).$$
\end{conj}

In the later discussion, for our own convenience, we will freely make linear changes of the variables for $\z^{SL_n}(s)$ and denote the resulting functions by $\z^{SL_n}(s)$ as well.

\subsection{$SL_2$}
By definition, a direct calculation (with a linear change of variable) leads to,
$$\z^{SL_2}(s):=\frac{\z(2s)}{1-q^{-2s+2}}+\frac{\z(2s-1)}{1-q^{2s}}.$$ Set $t=q^{-s}, T=t^2$ and $a_1=q+1-N$. 
Then 
$$\begin{aligned}\z^{SL_2}(t)=&
\frac{1-a_1T+qT^2}{(1-T)(1-qT)(1-q^2T)}+\frac{1-a_1qT+q^3T^2}{(1-qT)(1-q^2T^2)(1-\frac{1}{T})}\\
=&\frac{1-a_1T+qT^2-T(1-a_1qT+q^3T^2)}{(1-T)(1-qT)(1-q^2T)}.\end{aligned}$$
Consequently,
$$\z^{SL_2}(s)=\frac{1+(N-2)T+q^2T^2}{(1-T)(1-q^2T)}.$$

\begin{thm} Conjecture 7 holds for $SL_2$. That is to say,
the Uniformity and the Riemann Hypothesis hold for $\z^{SL_2}(s)$ and $\widehat\zeta_{E,2}(s)$.
\end{thm}

\noindent
{\it Remark.} Yoshida ([Y])  shows that, more generally, $\widehat\zeta_X^{SL_2}(s)$ satisfies the RH for all regular, integral, projective curve $X$ defined over $\mathbb F_q$.

\subsection{$SL_3$}
By definition, a direct calculation (with a linear change of variable) shows that
the zeta of $E$ associated to $SL_3$ is given by
$$\begin{aligned}\z^{SL_3}(s)=&\z(2)\cdot\Big(\frac{\z(3s)}{1-q^{-3s+3}}+\frac{\z(3s-2)}{1-q^{3s}}\Big)\\
&+\frac{\z(1)}{1-q^2}\cdot\Big(\frac{\z(3s)}{1-q^{-3s+2}}+\frac{\z(3s-2)}{1-q^{3s-1}}\Big)\\
&+\z(1)\cdot \frac{\z(3s-1)}{(1-q^{3s})(1-q^{-3s+3})}.\end{aligned}$$
With $t=q^{-s},\ T=t^3$, we have
$$\z^{SL_3}(s)=\frac{P_E^{SL_3}(T)}{(1-T)(1-qT)(1-q^2T)(1-q^3T)}$$ with $$\begin{aligned}&P_E^{SL_3}(T):=\Big(1+\frac{qN}{(q-1)(q^2-1)}\Big)\\
&\times\Big[(1+q^6T^4-(q^2+q+2-N)(T+q^3T^3)+2qT^2\big(1+(q+1-N)q\big)\Big]\\
&\qquad-\frac{N}{(q-1)(q^2-1)}\\
&\times\Big[1+q^6T^4-(q^3+2q+1-N)(T+q^3T^3)+2qT^2\big(1+(q^2(q+1-N)\big)\Big]\\
&\qquad-\frac{N}{q-1}\Big[(T+q^3T^3)-(q+1-N)qT^2\Big].\end{aligned}$$

\begin{lem} (i) There exists  a 
degree 2 polynomial $P_{E,o}^{SL_3}(T)$ of $T$ such that
$$P_E^{SL_3}(T)=(1-qT)(1-q^2T)P_{E,o}^{SL_3}(T);$$

\noindent
(ii) $P_{E,o}^{SL_3}(T)$ is given by
$$P_{E,o}^{SL_3}(T)=\Big(1+\frac{N}{q^2-1}\Big)+\Big[-2+\frac{2q-3}{q-1}N+\frac{N^2}{q^2-1}\Big]\cdot T+\Big(1+\frac{N}{q^2-1}\Big)q^3T^2.$$
\end{lem}

\noindent
{\it Proof.} (i) By functional equation, it suffices to show that $$P_E^{SL_3}(\frac{1}{q})=0.$$ Then routine checking.
In fact,  one can first set $qT=1$ in the above expression of $P_E^{SL_3}(T)$. Then using SIMPLIFY command of Mathematica, to verify that the resulting complicated combination in terms of $q$ and $N$ gives us 0 as wanted.

\noindent
(ii) You can directly calculate it by hands. Instead, first we have 
$$\begin{aligned}&P_E^{SL_3}(T):=\Big[1+\frac{N}{q^2-1}\Big]\cdot\Big(1+q^6T^4\Big)\\
&\times\Big[-(q^2+q+2)+N\frac{q-3}{q-1}+\frac{N^2}{q^2-1}\Big]\cdot\Big(T+q^3T^3\Big)\\
&\Big[2(q^2+q+1)-N\frac{2q^3-q^2-4q-3}{q^2-1}-\frac{N^2}{q-1}\Big]\cdot qT^2.\end{aligned}$$ Then using the PolynomialQuotientRemainder command of Mathematica, to divide $P_E^{SL_3}(T)$ by $q^3T^2-(q^2+q)T+1$. 
As a result, Mathematica would give us 
\vskip 0.30cm
\noindent
(a) the quotient, a degree two polynomial of $T$ with coefficients in terms of $q$ and $N$. Using Simplify commend of Matematica to get the result in the lemma;

\noindent
(b) the reminder, a linear polynomial  in  $T$ with coefficients in terms of $q$ and $N$. Using Simplify commend of Matematica to see that the coefficients of $T$ and constant terms are all 0.

This then completes the proof of the lemma.

\begin{thm}
(i) ({\bf Uniformity}$_3$) $$\widehat \zeta_{E,3}(s)=\z(1)\cdot\z^{SL_3}(s).$$
\noindent
(ii) ({\bf Riemann Hypothesis}$_3$) $$\widehat \zeta_{E,3}(s)=0\qquad\Rightarrow\qquad \mathrm{Re}(s)=\frac{1}{2}.$$
\end{thm}

\noindent
{\it Proof.} (i) is a direct consequence of  the closed formulas for $\widehat \zeta_{E,2}(s)$ and $\z^{SL_3}(s)$ in Prop 3 and Lem 9.

For (ii), it suffices to show that the discriminant of the degree two polynomial $P_{E,o}^{SL_3}(T)$ is strictly negative. Clearly,
$$\begin{aligned}\Delta_3=&\Big[-2+\frac{2q-3}{q-1}N+\frac{N^2}{q^2-1}\Big]^2-4\Big(1+\frac{N}{q^2-1}\Big)^2q^3\\
=&\Big[\Big(-2+\frac{2q-3}{q-1}N+\frac{N^2}{q^2-1}\Big)+2\Big(1+\frac{N}{q^2-1}\Big)q\sqrt q\Big]\\
&\times \Big[\Big(-2+\frac{2q-3}{q-1}N+\frac{N^2}{q^2-1}\Big)-2\Big(1+\frac{N}{q^2-1}\Big)q\sqrt q\Big].\end{aligned}$$ The first factor is strictly positive, while by Hasse's theorem for Artin zetas, the second factor is strictly negative. This then completes the proof.

\subsection{Distribution of zeros}

For  Artin zetas of elliptic curves $E/\mathbb F_q$, set
$$\cos\theta_{p}=\frac{p+1-N(E/\mathbb F_p)}{2\sqrt p},\qquad 0<\theta_p<\pi.$$

\begin{conj} ({\bf Sato-Tate Conjecture}) If  $\{E/\mathbb F_p\}$ are not (resulting from a global one of CM type), then
for $0\leq\alpha<\beta\leq \pi$,
$$\lim_{x\to \infty}\frac{\#\Big\{p\ \mathrm{prime}: p\leq x, \alpha\leq\theta_p\leq\beta\Big\}}{\#\big\{p\ \mathrm{prime}: p\leq x\big\}}=\frac{2}{\pi}\int_\alpha^\beta\sin^2 \theta d\theta.$$
\end{conj}

There are some exciting developments in this direction due to Taylor,  Clozel-Harris-Shepherd-Barron and Barnet-Lamb-Geraghty-Harris.

 Motivated by this, by the RH for $SL_2$ and $SL_3$, we set
 $$\cos\theta_{2,p}=\frac{N(E/\mathbb F_p)-2}{2p},$$
$$\cos\theta_{3,p}=\frac{-2+\frac{2p-3}{p-1}N(E/\mathbb F_p)+\frac{N(E/\mathbb F_p)^2}{p^2-1}}{2p\sqrt p\cdot \Big(1+\frac{N(E/\mathbb F_p)}{p^2-1}\Big)},$$
and hence to understand
$$\lim_{x\to \infty}\frac{\#\Big\{p\ \mathrm{prime}: p\leq x, \alpha\leq\theta_{2,p}\leq\beta\Big\}}{\#\big\{p\ \mathrm{prime}: p\leq x\big\}}$$ and 
$$\lim_{x\to \infty}\frac{\#\Big\{p\ \mathrm{prime}: p\leq x, \alpha\leq\theta_{3,p}\leq\beta\Big\}}{\#\big\{p\ \mathrm{prime}: p\leq x\big\}}.$$

By Hasse's theorem, we have $|N(E/\mathbb F_p)-p-1|\leq 2\sqrt p$. Thus,
$$\lim_{p\to\infty}\frac{N(E/\mathbb F_p)-2}{2p}=\frac{1}{2}$$ and
$$\lim_{p\to\infty}\frac{-2+\frac{2p-3}{p-1}N(E/\mathbb F_p)+\frac{N(E/\mathbb F_p)^2}{p^2-1}}{2p\sqrt p\cdot \Big(1+\frac{N(E/\mathbb F_p)}{p^2-1}\Big)}=0.$$
Consequently, $$\lim_{p\to\infty}\theta_{2,p}=\frac{\pi}{3},\qquad \lim_{p\to\infty}\theta_{2,p}=\frac{\pi}{2}.$$
Therefore, we have the following

\begin{prop} The distributions of zeros for rank 2 zeta, resp. rank 3 zetas of elliptic curves are of Dirac type. More precisely,
$$\lim_{x\to \infty}\frac{\#\Big\{p\ \mathrm{prime}: p\leq x, \alpha\leq\theta_{2,p}\leq\beta\Big\}}{\#\big\{p\ \mathrm{prime}: p\leq x\big\}}=\int_\alpha^\beta\delta_{\frac{\pi}{3}}\,dt$$ and 
$$\lim_{x\to \infty}\frac{\#\Big\{p\ \mathrm{prime}: p\leq x, \alpha\leq\theta_{3,p}\leq\beta\Big\}}{\#\big\{p\ \mathrm{prime}: p\leq x\big\}}=\int_\alpha^\beta\delta_{\frac{\pi}{2}}\,dt,$$
where $\delta_a$ denotes the Dirac distribution at $a$.
\end{prop}

\section{Counting Bundles}
Recall that the rank $r$ pure non-abelian zeta function of an elliptic curve $E/\mathbb F_q$ is given by
$$\widehat Z_{E,r}(t)
=\alpha_{E,r}\big(0\big)+\beta_{E,r}(0)\cdot \frac{(Q-1)T}{(1-T)(1-QT)}=\frac{P_{E,r}(T)}{(1-T)(1-QT)}$$
with $$P_{E,r}(T)=\alpha_{E,r}\big(0\big)-\Big[(Q+1)\alpha_{E,r}\big(0\big)-(Q-1)\beta_{E,r}\big(0\big)\Big]T+\alpha_{E,r}\big(0\big) QT^2.$$
Here $t=q^{-s},\, Q=q^r$ and $T=t^r$. Thus to determine it, we need to know the invariants 
$\alpha_{E,r}(0)$ and $\beta_{E,r}(0)$.

As said, the $\beta$-invariant has been studied by many authors. In fact Harder-Narasimhan, Desale-Ramanan, and Zagier's formula can be arranged as follows:

\begin{thm} $$\beta_{E,n}(0)=\sum_{n_1+\dots+n_k=n}\prod_{j=1}^{k-1}\frac{1}{q^{n_j+n_{j+1}}-1}
v_{n_1,\dots,n_k}$$ where
$$v_{n_1,\dots,n_k}:=\prod_{j=1}^kv_{n_j}\qquad \mathrm{with}\qquad v_n:=\widehat\zeta_E(1)\widehat\zeta_E(2)\cdots \widehat\zeta_E(n).$$
\end{thm}

This is significantly clearer than the original formula stated in Theorem 1, since the parabolic reduction
structure appears and the parabolic coefficients $$e_{n_1,\dots,n_k}=\prod_{j=1}^{k-1}\frac{1}{q^{n_j+n_{j+1}}-1}$$ is {\it environmentally free}, i.e., only determined by the group structure but independent of curves. Indeed, for general genus curve $X$ the same formula holds if we rewrite the left hand as $$u_n:=\frac{\beta_{X,n}(0)}{q^{\frac{n(n+1)}{2}(g-1)}}.$$

\begin{thm} (See e.g., [W5])
$$u_n=\sum_{n_1+\dots+n_k=n}(-1)^ke_{n_1,\dots,n_k}\cdot v_{n_1,\dots,n_k}.$$
\end{thm}

As for the $\alpha$-invariant, by our previous calculation in lower ranks, we introduce the
following

\begin{conj} ({\bf Counting Miracle}) For elliptic curves $E/\mathbb F_q$,
$$\alpha_{E,n+1}(0)=\beta_{E,n}(0)$$
\end{conj}
We have checked it  for $n=1,2,3,4,5$. 

To understand this, let us introduce the so-called Atiyah bundles $I_r$ inductively.
The starting point is $I_1=\mathcal O_E$. Then we consider the extension $\mathcal O_E$ by $\mathcal O_E$.
Since $\mathrm{Ext}_E^1(\mathcal O_E,\mathcal O_E)$ is 1 dimensional,  there is, up to isomorphism, only one non-trivial extension. This is $I_2$, namely, we have the non-trivial extension
$$0\to \mathcal O_E\to I_2\to\mathcal O_E\to 0.$$ Inductively, we know that $\mathrm{Ext}_E^1(I_{r-1},\mathcal O_E)$ is 1 dimensional, so there is, up to isomorphism, only one non-trivial extension of $I_{r-1}$ by $\mathcal O_E$. This is $I_r$, namely, we have the non-trivial extension
$$0\to \mathcal O_E\to I_r\to I_{r-1}\to 0.$$

\begin{lem} (i) ([A, Thm 8]) $h^0(E,I_r)=1$ and 
$$I_r\otimes I_s=I_{r-s+1}\oplus I_{r-s+3}\oplus\cdots\oplus I_{r+s},\qquad r\geq s;$$

\noindent
(ii) For a partition $n=m_1\cdot r_1+m_2\cdot r_2+\cdots+m_s\cdot r_s$
$$=
(r_1+\cdots+r_1)+(r_2+\cdots+r_2)+\cdots+(r_s+\cdots+r_s)$$ arranging in the order 
$r_1<r_2<\dots<r_s$, we have
$$\begin{aligned}\#\mathrm{Aut}\Big(\oplus_{j=1}^sI_{r_j}^{\oplus m_j}\Big)&=q^{2\sum_{1\leq i<j\leq s}r_im_im_j}\\
&\times\prod_{j=1}^s(q^{m_j}-1)(q^{m_j}-q)\cdots(q^{m_j}-q^{m_j-1})q^{m_j^2(r_j-1)}.\end{aligned}$$
\end{lem}

\noindent
{\it Proof.} In (i), by definition, the first on $h^0$ is obvious, and the multiplicative relation is given in  Thm 8 of [A]. As for (ii), we use the natural surjective morphism $$\mathrm{Aut}\Big(\oplus_{j=1}^sI_{r_j}^{\oplus m_j}\Big)
\buildrel\Pi\over\to\prod_{j=1}^s \mathrm{Aut}\Big(I_{r_j}^{\oplus m_j}\Big)$$ with kernel
$$\Big[\mathrm{Id}+\oplus_{i<j}\mathrm{Hom}\Big(I_{r_i}^{\oplus m_i},I_{r_j}^{\oplus m_j}\Big)\Big]
\times \Big[\mathrm{Id}+\oplus_{i<j}\mathrm{Hom}\Big(I_{r_j}^{\oplus m_j},I_{r_i}^{\oplus m_i}\Big)\Big]$$

Note that $$\mathrm{Aut}\Big(I_{r}\Big)\simeq\Big\{
\begin{pmatrix}a&b_1&b_2&\cdots&b_{r-2}&b_{r-1}\\
0&a&b_1&\cdots&b_{r-3}&b_{r-2}\\
0&0&a&\cdots&b_{r-4}&b_{r-3}\\
\cdots&\cdots&\cdots&\cdots&\cdots&\cdots\\
0&0&0&\cdots&a&b_1\\
0&0&0&\cdots&0&a\end{pmatrix}\in GL_r(\mathbb F_q)\Big\}.$$
Hence $$\#\mathrm{Aut}\Big(I_{r}\Big)=(q-1)q^{r-1}.$$
More generally, for $\mathrm{Aut}\Big(I_{r_j}^{\oplus m_j}\Big)$, let us decompose its elements
into $m_j\times m_j$ blocks of size $r_j\times r_j$. Then corresponding to the simple factor $q-1$ for $I_r$, now we have the factor $(q^{m_j}-1)(q^{m_j}-q^2)\cdots(q^{m_j}-q^{m_j-1})$, the number elements of $GL_{m_j}(\mathbb F_q)$, and the offer diagonal parts of $I_r$ now give us the total number
$(q^{r_j-1})^{m_j^2}$ since there are $m_j\times m_j$-blocks.
Consequently, we have
 $$\#\mathrm{Aut}\Big(I_{r_j}^{\oplus m_j}\Big)=\Big[(q^{m_j}-1)(q^{m_j}-q)\cdots(q^{m_j}-q^{m_j-1})\Big]\cdot q^{m_j^2(r_j-1)}.$$ 
 Thus to complete the proof, we need to how that
 $$\#\oplus_{i<j}\mathrm{Hom}\Big(I_{r_i}^{\oplus m_i},I_{r_j}^{\oplus m_j}\Big)=
 q^{\sum_{1\leq i<j\leq s}r_im_im_j}.$$
This results from Atiyah's multiplicative structure on $I_r$.
Indeed, 
$$\begin{aligned}\mathrm{Hom}\Big(I_{r_i}^{\oplus m_i},I_{r_j}^{\oplus m_j}\Big)\simeq& H^0(E,I_{r_i}^\vee\otimes I_{r_j})^{\oplus m_im_j}\\
\simeq& H^0(E,I_{r_i}\otimes I_{r_j})^{\oplus m_im_j}\\
\simeq& H^0(E,I_{r_j-r_i+1}\oplus I_{r_j-r_i+3}\oplus\cdots\oplus I_{r_j+r_i-1})^{\oplus m_im_j}\\
=&(\mathbb F_q)^{r_im_im_j}.\end{aligned}$$
This then completes the proof.

With this in mind, now we introduce the following

\begin{conj} ({\bf Miracle of $I_r$})
$$\begin{aligned}\sum_{{\sum_{i=1}^sr_im_i=n+1}}&\frac{q^{h^0(E,\oplus_{j=1}^sI_{r_j}^{\oplus m_j})}-1}{\#\mathrm{Aut}\Big(\oplus_{j=1}^sI_{r_j}^{\oplus m_j}\Big)}
=\sum_{{\sum_{i=1}^sr_im_i=n}}\frac{1}{\#\mathrm{Aut}\Big(\oplus_{j=1}^sI_{r_j}^{\oplus m_j}\Big)}\\
=&\frac{q^{\frac{n(n+1)}{2}}}{(q^{n+1}-1)(q^n-1)(q^{n-1}-1)\cdots(q-1)}.\end{aligned}$$
\end{conj}

\begin{prop} Miracle of $I_r$ implies the counting miracle.
\end{prop}

\noindent
{\it Proof.} Let $V$ be a semi-stable vector bundle of rank $r$ over $E/\mathbb F_q$.
Then the graded bundle associated to its Jordan-H\"older filtrations decomposes as
$$G(V)=L_1\oplus L_2\oplus\cdots\oplus L_r$$ with $L_i$'s line bundles of  degree on $\overline E:=E\otimes_{\mathbb F_q}\overline{\mathbb F_q}$. As $L_i$ need not be defined over $\mathbb F_q$, usually it is a bit complicated to classify $V$. As a matter of fact, this classification problem is related with arithmetic of elliptic curve, say, depending on the number of $\mathbb F_q$-rational $r$ torsions of $E$.
Instead of counting them using a complete list, we first note that in order to have a non-trivial contribution to
$\alpha$, $h^0(V)\not=0$. Guided by this, we regroup the summation as follows:
$$\alpha_{E,r}(0)=\sum_{i=1}^r\sum_{\substack{V,\\ G(V)=\mathcal O_E^{\oplus i}\oplus L_{i+1}\oplus \cdots\oplus L_r,\\ L_j\not=\mathcal O_E}}\frac{q^{h^0(E,V)}-1}{\#\mathrm{Aut}V}.$$
With $G(V)=\mathcal O_E^{\oplus i}\oplus L_{i+1}\oplus \cdots\oplus L_r, L_j\not=\mathcal O_E$, since there is no morphism between $\mathcal O_E$ and $L_j$'s,
$$\mathrm{Aut}V\simeq\mathrm{Aut}\,U\times \mathrm{Aut}W$$ with $G(U)=\mathcal O_E^{\oplus i}$ and $G(W)=L_{i+1}\oplus\cdots\oplus L_r$, $L_j\not=\mathcal O_E$.
Consequently, $$\begin{aligned}\alpha_{E,r}(0)=&\sum_{i=1}^r\sum_{\substack{V,\, V=U\oplus W,\\ G(U)=\mathcal O_E^{\oplus i},\\ G(W)= L_{i+1}\oplus \cdots\oplus L_r,\\ L_j\not=\mathcal O_E}}\frac{q^{h^0(E,U)}-1}{\#\mathrm{Aut}\,U\cdot \#\mathrm{Aut}W}\\
=&\sum_{i=1}^r
\sum_{\substack{U,\\G(U)=\mathcal O_E^{\oplus i}}}\frac{q^{h^0(E,U)}-1}{\#\mathrm{Aut}\,U}
\sum_{\substack{W,\\ G(W)= L_{i+1}\oplus \cdots\oplus L_r,\\ L_j\not=\mathcal O_E}}\frac{1}{\#\mathrm{Aut}W}.\end{aligned}$$
Now assume that we have the Miracle for $I_r$, then 
$$\sum_{\substack{U,\\G(U)=\mathcal O_E^{\oplus i}}}\frac{q^{h^0(E,U)}-1}{\#\mathrm{Aut}\,U}
=\sum_{\substack{U,\\G(U)=\mathcal O_E^{\oplus i-1}}}\frac{1}{\#\mathrm{Aut}\,U}.$$
Therefore, by reversing the above discussion with $q^{h^0}-1$ replaced by $1$, we get
$$\begin{aligned}\alpha_{E,r}(0)
=&\sum_{i=1}^r
\sum_{\substack{U,\\G(U)=\mathcal O_E^{\oplus i-1}}}\frac{1}{\#\mathrm{Aut}\,U}
\sum_{\substack{W,\\ G(W)= L_{i+1}\oplus \cdots\oplus L_r,\\ L_j\not=\mathcal O_E}}\frac{1}{\#\mathrm{Aut}W}\\
=&\sum_{i=1}^r\sum_{\substack{V,\, V=U\oplus W,\\ G(U)=\mathcal O_E^{\oplus i-1},\\ G(W)= L_{i+1}\oplus \cdots\oplus L_r,\\ L_j\not=\mathcal O_E}}\frac{1}{\#\mathrm{Aut}\,U\cdot \#\mathrm{Aut}W}\\
=&\sum_{i=1}^r\sum_{\substack{V,\\ G(V)=\mathcal O_E^{\oplus i-1}\oplus L_{i+1}\oplus \cdots\oplus L_r,\\ L_j\not=\mathcal O_E}}\frac{1}{\#\mathrm{Aut}V}\\
=&\beta_{E,r-1}(0).\end{aligned}$$

\section{Distribution of Zeros}
\subsection{The Riemann Hypothesis}
For high rank pure zeta functions of elliptic curves, with our works on rank two and three zeta functions, we now introduce the following
\begin{conj} ({\bf Riemann Hypothesis}) 
$$\widehat \zeta_{E,r}(s)=0\qquad\Rightarrow\qquad \mathrm{Re}(s)=\frac{1}{2}.$$
\end{conj}

Since $$\begin{aligned}\widehat Z_{E,r}(t)
=&\alpha_{E,r}\big(0\big)+\beta_{E,r}(0)\cdot \frac{(Q-1)T}{(1-T)(1-QT)}\\
=&\frac{P_{E,r}(T)}{(1-T)(1-QT)}\end{aligned}$$
with $$\begin{aligned}P_{E,r}(T)=&\alpha_{E,r}\big(0\big)+\alpha_{E,r}\big(0\big) QT^2\\
&-\Big[(Q+1)\alpha_{E,r}\big(0\big)-(Q-1)\beta_{E,r}\big(0\big)\Big]T.\end{aligned}$$
For our own use, set
$$a_{E,r}=\frac{\beta_{E,r}\big(0\big)}{\alpha_{E,r}\big(0\big)}.$$

The Riemann Hypothesis means that
$$\Delta_r=\Big[(Q+1)\alpha_{E,r}\big(0\big)-(Q-1)\beta_{E,r}\big(0\big)\Big]^2-4\alpha_{E,r}\big(0\big)^2Q<0,$$ or the same,
$$\begin{aligned}0>&\Big[(Q+1)\alpha_{E,r}\big(0\big)-(Q-1)\beta_{E,r}\big(0\big)+2\alpha_{E,r}\big(0\big)\sqrt Q\Big]\\
&\cdot \Big[(Q+1)\alpha_{E,r}\big(0\big)-(Q-1)\beta_{E,r}\big(0\big)-2\alpha_{E,r}\big(0\big)\sqrt Q\Big]\\
=&(Q-1)\alpha_{E,r}\big(0\big)^2
\Big[\sqrt Q\Big(1-a_{E,r}\Big)+\Big(1+a_{E,r}\Big)
\Big]\cdot \Big[\sqrt Q\Big(1-a_{E,r}\Big)-\Big(1+a_{E,r}\Big)
\Big].\end{aligned}
$$ 
Or, equivalently,
$$\begin{cases} \sqrt Q\big(1-a_{E,r}\big)+\big(1+a_{E,r}\big)>0\\
\sqrt Q\big(1-a_{E,r}\big)-\big(1+a_{E,r}\big)<0.\end{cases}$$
That is to say,
$$\frac{\sqrt Q-1}{\sqrt Q+1}<a_{E,r}< \frac{\sqrt Q+1}{\sqrt Q-1}.$$

\begin{conj} ({\bf Riemann Hypothesis}) $$1-\frac{2}{\sqrt {q^r}+1}<a_{E,r}< 1+\frac{2}{\sqrt {q^r}-1}.$$
 In particular,
$$a_{E,r}=\frac{\beta_{E,r}\big(0\big)}{\beta_{E,r-1}\big(0\big)}\to 1,\qquad q\cdot r\to\infty.$$
\end{conj}

Discussion: For $r\to \infty$, we expect that
$$\beta_{E,r}(0)\,\sim\, \widehat\zeta_E(1)\widehat\zeta_E(2)\cdots \widehat\zeta_E(r).$$ 
So $$a_{E,r}\,\sim \,\widehat\zeta_E(r)=1+\frac{Nq^{r-1}}{(q^r-1)(q^{r-1}-1)}$$

This implies asymptotically the RH holds.

\subsection{Distribution of Zeros}
To end this paper, let us consider the distribution of zeros under the assumption that the RH holds. Set
$$\cos\theta_{r,p}^E:=\frac{(p^r-1)\beta_{E,r}\big(0\big)+(p^r+1)\alpha_{E,r}\big(0\big)}{2\sqrt{p^r}\alpha_{E,r}\big(0\big)}$$
Then $$\cos\theta_{r,p}^E
=\frac{(p^r-1)a_{E,r}+(p^r+1)}{2\sqrt{p^r}}
=\frac{p^r(1-a_{E,r})+(1+a_{E,r})}{2\sqrt{p^r}}.$$ Thus
in assuming the Riemann Hypothesis, we have
$$\lim_{pr\to\infty}\cos\theta_{r,p}^E= 0.$$ Therefore,
$$\lim_{x\to \infty}\frac{\#\Big\{p\ \mathrm{prime}: p\leq x, \alpha\leq\theta_{r,p}\leq\beta\Big\}}{\#\big\{p\ \mathrm{prime}: p\leq x\big\}}=\int_\alpha^\beta\delta_{\frac{\pi}{2}}\,dt,$$
where $\delta_a$ denotes the Dirac distribution at $a$.
So the distribution of zeros of high rank zetas are very much different from that of Artin's.


As we expect that, asymptotically, $$a_{E,r}\sim\widehat\zeta_E(r)=1+\frac{Np^{r-1}}{(p^r-1)(p^{r-1}-1)},$$ for the distributions of zeros of high rank zetas, with $\delta_{\frac{\pi}{2}}$ understood, 
we should go further to analysis the subdominant term, in order to see the refine structure of the zeros.


\vskip 2.0cm
\centerline{\bf REFERENCES}
\vskip 0.30cm
\noindent
[A] M.F. Atiyah, Vector bundles over an elliptic curve, Proc. London Math Soc (3) 7 (1957), 414-452 
\vskip 0.10cm
\noindent
[DR] U.V. Desale \& S. Ramanan, Poincare polynomials of the variety of stable bundles, Math. Ann 26 (1975) 233-244
\vskip 0.10cm
\noindent[HN] G. Harder \& M.S. Narasimhan, On the cohomology groups of moduli spaces of vector bundles on curves, Math. Ann. 212, 215-248 (1975)
\vskip 0.10cm
\noindent
[W1] L. Weng, Non-abelian zeta function for function fields, Amer. J. Math., 127 (2005), 973-1017
\vskip 0.10cm
\noindent
[W2] L. Weng, A geometric approach to $L$-functions, in {\it Conference on L-Functions}, pp. 219-370, World Sci (2007)
\vskip 0.10cm
\noindent
[W3] L. Weng, Symmetry and the Riemann Hypothesis, in {\it Algebraic and Arithmetic Structures of Moduli Spaces}, ASPM 58, 173-223 (2010)
\vskip 0.10cm
\noindent
[W4] L. Weng, Counting Bundles, preprint, 2011 
\vskip 0.10cm
\noindent
[W5] L. Weng, Parabolic Reduction, Stability and the Mass, in preparation
\vskip 0.10cm
\noindent
[Y] H. Yoshida, manuscripts, Dec., 2011
\vskip 0.10cm
\noindent
[Z] D. Zagier, Elementary aspects of the Verlinde formula and the Harder-Narasimhan-Atiyah-Bott formula, in {\it Proceedings of the Hirzebruch 65 Conference on Algebraic Geometry}, 445-462 (1996)
\vskip 3.5cm
\noindent
{\bf Lin WENG}\footnote{{\bf Acknowledgement.}  We would like to thank H. Yoshida for his keen interests in our works. 
His related works on zetas and their zeros is one of our starting points for this new attempt to understand zetas associated to function fields.

This work is partially supported by JSPS.}


\noindent
Graduate School of Mathematics, Kyushu University,
 Fukuoka 819-0395

\noindent
E-Mail: weng@math.kyushu-u.ac.jp
\end{document}